\documentclass[12pt,reqno]{amsart}

\usepackage{amssymb,amsmath}
\usepackage{graphicx}

\newtheorem{lemma}{Lemma}[section]

\theoremstyle{definition}

\newtheorem{example}[lemma]{Example}

\theoremstyle{remark}

\newcommand{\Sing}{\operatorname{Sing}}

\title{Polar varieties revisited}

\author{Ragni Piene}
\address{Department of Mathematics, University of Oslo,\\
P.O.Box 1053 Blindern, NO-0316 Oslo, Norway}
\email{ragnip@math.uio.no}

\begin{document}

\begin{abstract}
We recall the definition of classical polar varieties, as well as those of affine and projective reciprocal polar varieties. The latter are defined with respect to a non-degenerate quadric, which gives us a notion of orthogonality. In particular we relate the reciprocal polar varieties to the  ``Euclidean geometry'' in projective space. The Euclidean distance degree and the degree of the focal loci can be expressed in terms of the ranks, i.e., the degrees of the classical polar varieties, and hence these characters can be found also for singular varieties, when one can express the ranks in terms of the singularities. 
\end{abstract}

\maketitle

\section{Introduction}

The theory of polars and polar varieties has played an 
important role in the quest for understanding and classifying projective
varieties. Their use in the definition of projective invariants is the
very
basis for the geometric approach to the theory of characteristic classes, such
as Todd classes and Chern classes. In particular this approach gives a way of defining Chern classes for \emph{singular} projective varieties (see e.g. \cite{Piene1,Piene3}).
The \emph{local} polar varieties were used by L\^e--Teissier, Merle and others in the  study of singularities.

More recently, polar varieties have been applied to study
 the topology of real affine varieties and to find real
solutions of polynomial equations 
 by Bank et al., Safey El Din--Schost, and Mork--Piene
 \cite{Bank,Bank2,Bank3,Bank4,Safey,Mork,MorkPiene},
to complexity questions by B\" urgisser--Lotz \cite{BurgL},
to foliations by Soares \cite{Soares} and others, to
 focal loci and caustics of reflection by Catanese and Trifogli \cite{C-T,T} and
 Josse--P\`ene \cite{JP}, to
 Euclidean distance degree by Draisma et al. \cite{ED}.

In this note I will explore the relation of polar and reciprocal polar varieties  of possibly singular varieties to   the Euclidean normal bundle, the Euclidean distance degree, and the focal loci. For simplicity, we work over an algebraically closed field of characteristic $0$.

\section{Classical polar varieties}

Let $\mathbb G(m,n)$ denote the  Grassmann variety of  $m$-dimensional linear subspaces of $\mathbb P^n$.
Let $L_k\subset \mathbb P^n$ be a linear subspace of dimension $n-m+k-2$.
 Consider the special Schubert variety
\[\Sigma(L_k):=\{W\in \mathbb G(m,n) \,|\, \dim W\cap L_{k}\ge k-1\}.\]
It has a natural structure as a determinantal scheme and is of codimension $k$ in $\mathbb G(m,n)$.

 As is well known,  if $L_k'$ is another linear subspace of dimension $n-m+k-2$, then $\Sigma(L_k')$ and  $\Sigma(L_k)$ are projectively equivalent (in particular, their rational equivalence classes are equal). So we will often just write $\Sigma_k$ for such a variety.

\begin{example}
The case $m=1$, $n=3$. For $\mathbb G(1,3)=\{\text{lines in } \mathbb P^3\}$, the special Schubert varieties are
\begin{align*}
\Sigma_1&=\{\text{lines meeting a given line}\} \\
\Sigma_2&=\{\text{lines in  a given plane}\}.
\end{align*}
\end{example}
\begin{example}
The case $m=2$, $n=5$. 
For  $\mathbb G(2,5)=\{\text{planes in } \mathbb P^5\}$, the special Schubert varieties are
\begin{align*}
\Sigma_1&=\{\text{planes meeting  a given plane}\} \\
\Sigma_2&=\{\text{planes intersecting a given 3-space in a line}\}\\
\Sigma_3&=\{\text{planes contained in a given hyperplane}\}
\end{align*}
\end{example}

More general Schubert varieties are defined similarly, by giving conditions with respect to flags of linear subspaces (see e.g. \cite{KL}). For example, in Example 1, we could consider 
\[\Sigma_{0,2}=\{\text{lines in a given plane through a given point in the plane}\}.\]

Let $X\subset \mathbb P^n$ be a (possibly singular) variety of dimension $m$. The \emph{Gauss map} $\gamma\colon X\dashrightarrow \mathbb G(m,n)$ is the rational map that sends a nonsingular point $P\in X_{\rm sm}$ to the (projective) tangent space $T_PX$, considered as a point in the Grassmann variety. More precisely, if $V=H^0(\mathbb P^n,\mathcal O_{\mathbb P^n}(1))$, then $\gamma$ is given on $X_{\rm sm}:=X\setminus \Sing X$ by the restriction of the quotient
\[V_X\to \mathcal P_X^1(\mathcal L),\] where $\mathcal P_X^1(\mathcal L)$ denotes the bundle of 1st principal parts of the line bundle $\mathcal L:=\mathcal O_{\mathbb P^n}(1)|_X$. Note that restricted to $X_{\rm sm}$ the sheaf of principal parts is locally free with rank $m+1$, and that the fibers (over $X_{\rm sm}$) of $\mathbb P(\mathcal P_X^1(\mathcal L)) \subset X\times \mathbb P^n$ over $X$ (with respect to projection on the first factor) define the (projective) tangent spaces of $X$.

\medskip

The \emph{polar varieties} of $X\subset \mathbb P^n$ are the closures of the inverse images
\[M_k:=\overline{\gamma|_{X_{\rm sm}}^{-1}\Sigma_k}\]
of the special Schubert varieties via the Gauss map \cite[p.~252]{Piene1}.
\medskip

In the situation of Example 1,
 let $X\subset \mathbb P^3$ be a curve. Then $M_1$ is the set of nonsingular points $P\in X$ such that the tangent line $T_PX$ meets a given line, i.e., it is the ramification points of the linear projection map $X\to \mathbb P^1$.
\medskip

In Example 2, let $X\subset \mathbb P^5$ be a surface. Then $M_1$ is the ramification locus of the projection map $X\to \mathbb P^2$ with center a plane, and
$M_2$ consists of points $P\in X_{\rm sm}$ such that the tangent plane $T_PX$ intersects a given $3$-space in a line. One could also consider  more general polar  varieties, corresponding to general Schubert varieties, like  the set of points $P$ such that $T_PX$ meets a given line; this is the ramification locus of the linear projection of $X$ to $\mathbb P^3$ with the line as center. Note that the cardinality of this 0-dimensional variety is  equal to the degree of the tangent developable of $X$.

\section{Polar classes and Chern classes}
It was noted classically that the polar varieties are invariant under linear projections and sections. Therefore the \emph{polar classes}, i.e., their rational equivalence classes, are \emph{projective invariants} of the variety. Already Noether, Segre and Zeuthen observed that certain integral combinations of the polar classes of surfaces are \emph{intrinsic} invariants, i.e., they depend only on the surface, not on the given projective embedding. This was pursued by Todd and Severi, who used the polar classes to define what are now called the Chern classes. The formula is the following:
\[c_k(X)=\sum_{i=0}^k (-1)^{k-i}\binom{m+1-k+i}{i}[M_i]h^i,\]
where $h$ denotes the class of a hyperplane section. Since this expression makes sense also for singular varieties, it gives a definition of Chern classes for singular projective varieties,  called the \emph{Chern--Mather} classes (see e.g. \cite{Piene3}). The formula can be inverted to express the polar classes in terms of Chern classes. When the variety $X$ is nonsingular, this is just the expression coming from the fact that in this case,
\[[M_k]=c_k(\mathcal P_X^1(\mathcal L)).\]
Using the canonical exact sequence
\[0\to \Omega_X^1\otimes \mathcal L \to \mathcal P_X^1(\mathcal L) \to \mathcal L \to 0,\]
we compute, with $c_i(X)=c_i({\Omega_X^1}^\vee)=(-1)^ic_i(\Omega_X^1)$,
\[[M_k]=\sum_{i=0}^k (-1)^{k-i}\binom{m+1-k+i}ic_{k-i}(X)h^i.\]
When $X$ is singular, we replace $X$ by its Nash transform $\overline X$, i.e., $\nu\colon\overline X\to X$ is the smallest proper modification of $X$ such that $\gamma$ extends to a morphism $\overline \gamma\colon \overline X \to \mathbb G(m,n)$. If we denote by $V_{\overline X} \to\mathcal P$ the locally free quotient corresponding to $\overline \gamma$, then it follows from the definition that we have
\[ [M_k]=\overline \gamma_*c_k(\mathcal P).\]
In many cases, this allows us to find formulas for the degrees of the polar classes in terms of the singularities of $X$, see \cite{Piene1,Piene2,Piene3}.

\section{Reciprocal polar varieties}

In \cite{Bank3,Bank4} Bank et al. introduced what they called \emph{dual} polar varieties. These varieties were further studied in \cite{Mork,MorkPiene} under the name of \emph{reciprocal} polar varieties. They are defined with respect to a non degenerate quadric in the ambient projective space, and sometimes with respect to the choice of a hyperplane at infinity.

Let $Q\subset \mathbb P^n$ be a non degenerate quadric. Then $Q$ induces a \emph{polarity}, classically called a \emph{reciprocation}, between points and hyperplanes in $\mathbb P^n$. The \emph{polar} hyperplane $P^\perp$ of $P$ is the linear span of the points on $Q$ such that the tangent hyperplane to $Q$ at that point contains $P$:
\[P=(b_0:\cdots :b_n)\mapsto P^\perp \colon \sum_{i=0}^n b_i\frac{\partial q}{\partial X_i}=0.\]
The map $P \mapsto P^\perp$ is the isomorphism $\mathbb P^n\to ( \mathbb P^n)^\vee$ given by the symmetric bilinear form associated with the quadratic polynomial $q$ defining $Q$.  The point $P^\perp$ in the dual projective space represents a hyperplane in $\mathbb P^n$. If $L\subset \mathbb P^n$ is a linear space of dimension $r$, then  $L^\perp \subset (\mathbb P^n)^\vee$ gives a $(n-r-1)$-dimensional linear subspace $(L^\perp)^\vee \subset \mathbb P^n$, which we will, by abuse of notation, also denote by $L^\perp$.
 Note that if $H$ is a hyperplane, its \emph{pole} $H^\perp$ is the intersection of the tangent hyperplanes to $Q$ at the points of intersection with $H$. The map $H\mapsto H^\perp$ is the inverse of the above map $P\mapsto P^\perp$.

If the quadric is $q=\sum_{i=0}^n X_i^2$, then the polar of 
\[P=(b_0:\cdots :b_n)\] 
is the hyperplane 
\[P^\perp: b_0X_0+\ldots +b_nX_n=0.\]

Let $X\subset \mathbb P^n$ be a (possibly singular) variety of dimension $m$. Let
$K_i\subset \mathbb P^n$ be a linear subspace of dimension $i+n-m-1$. Recall \cite[Def. 2.1, p.~104]{MorkPiene}, \cite[p.~527]{Bank3} the
definition of the $i$-th reciprocal polar variety of $X$ (with respect to $K_i$):
 $W_{K_i}^{\perp}(X)$, $1\leq i \leq m$, 
is  the Zariski
closure of the set
\begin{displaymath}
\{ P\in X_{\rm sm}\setminus K_i^{\perp}\,|\, T_P X \not\pitchfork
\langle P,K_i^{\perp}\rangle ^{\perp}\},
\end{displaymath}
where $T_PX$ denotes the tangent space to $X$ at the point $P$ and
 the notation $M\not\pitchfork N$ means that the two linear spaces $M$, $N$ are not transversal, i.e., that their linear span $\langle M,N\rangle$ is not equal to the whole ambient space $\mathbb P^n$.
For general $K_i$, the $i$th reciprocal variety has codimension $i$.

The condition $T_P X \not\pitchfork \langle P,K_i^{\perp}\rangle ^{\perp}$ is equivalent to the condition $\dim (T_PX\cap \langle P,K_i^{\perp}\rangle ^{\perp})\ge i-1$, and $T_PX\cap \langle P,K_i^{\perp}\rangle ^{\perp}$ is equal to $T_PX\cap P^\perp\cap K_i$, or to $\langle T_PX^\perp,P\rangle ^\perp\cap K_i$, or to $\langle\langle T_PX^\perp,P\rangle,K_i^\perp \rangle ^\perp$. That the dimension of the latter space is greater or equal to $i-1$ is equivalent to 
$\dim \langle \langle T_PX^\perp,P\rangle,K_i^\perp \rangle \le n-i$, or to $\dim (\langle T_PX^\perp,P\rangle \cap K_i^\perp)\ge 0$.
Thus we  obtain the following description of the $i$th polar variety:
\[W_{K_i}^{\perp}(X)=\overline{\{P\in X_{\rm sm}\setminus K_i^{\perp}\,|\, \dim (\langle T_PX^\perp,P\rangle \cap K_i^\perp)\ge 0 \}}.\]
In the case $i=m$ we have that $K_m$ is a hyperplane, so $K_m^\perp$ is a point. Assuming this is not a point on $X$, we see that the $m$th reciprocal polar variety of $X$ with respect to $K_m$ is the (finite) set of nonsingular points $P\in X$ 
\[W_{K_m}^\perp(X)=\{P\in X_{sm} \, | \, K_m^\perp \in \langle T_PX^\perp,P\rangle \}.\] 
\medskip

Let $H_\infty =Z(x_0) \subset \mathbb P^n$ denote ``the hyperplane at infinity'', set $\mathbb A^n =\mathbb P^n \setminus H_\infty$, and consider the affine part $Y:=X\cap \mathbb A^n$ of $X$. Let $L=K_i\subseteq H_\infty$ be a linear space of dimension $i+n-m-1$. Define the \emph{affine reciprocal polar variety} to be the affine part of the reciprocal polar variety:
\[W^\perp_L(Y):=W^\perp_L(X)\cap \mathbb A^n.\]
 The linear variety $\langle P,L^\perp\rangle^\perp$ is contained in the hyperplane $H_\infty$, so we can consider the affine cone $I_{P,L^\perp}$ of $\langle P,L^\perp\rangle^\perp$ as a linear variety in the affine space $\mathbb A^n$. Then the affine reciprocal polar variety can be written as
\[W^\perp_L(Y)=\overline{\{P\in Y_{\rm sm}\setminus L^\perp \,|\, t_PY\not\pitchfork I_{P,L^\perp}\}},\]
where $t_PY$ denotes the affine tangent space at $P$, translated to the origin.

Assume $X=Z(F_1,\ldots,F_r)$ (for some $r\ge n-m$). Consider the case $L=K_m=H_\infty$. Then we see  that $W^\perp_{H_\infty}(Y)$ is the closure of the set of smooth points of $Y$ where the $(n-m+1)$-minors of the matrix
\[
\left(\begin{array}{ccc}
\frac{\partial q}{\partial x_1} &\cdots & \frac{\partial q}{\partial x_n}\\
\frac{\partial F_1}{\partial x_1} & \cdots & \frac{\partial F_1}{\partial x_n}\\
\vdots & \vdots &\vdots\\
\frac{\partial F_r}{\partial x_1} & \cdots & \frac{\partial F_r}{\partial x_n}.
\end{array}\right)
\]
vanish.
This generalizes \cite[Prop. 3.2.5]{MorkPiene}.

\medskip

\begin{example}
Assume $q$ is given by $q(x_0,x_1,\ldots,x_n)=x_0^2+\sum_{i=1}^n(x_i-a_ix_0)^2$ for some $a_1,\ldots,a_n$ such that $\sum_{i=1}^n a_i\neq 1$. Then $q$ restricts to (essentially) the square of the Euclidean distance function on $\mathbb A^n=\mathbb P^n \setminus H_\infty$, namely to $\sum_{i=1}^n(x_i-a_i)^2+1$. The affine reciprocal polar variety is given (on smooth points) by the vanishing of the $(n-m+1)$-minors of the matrix
\[
\left(\begin{array}{ccc}
x_1-a_1&\cdots & x_n-a_n\\
\frac{\partial F_1}{\partial x_1} & \cdots & \frac{\partial F_1}{\partial x_n}\\
\vdots & \vdots &\vdots\\
\frac{\partial F_r}{\partial x_1} & \cdots & \frac{\partial F_r}{\partial x_n}.
\end{array}\right)
\]
just as in \cite[(2.1)]{ED}.
\end{example}

\section{Euclidean normal bundles}
 Consider a variety $X\subset \mathbb P^n$ 
 and let $Q\subset \mathbb P^n$ be a non-degenerate quadric. As we saw in the previous section, the quadric induces a polarity on $\mathbb P^n$, which can be viewed as an orthogonality, like what one has in a Euclidean space. In \cite{C-T} this was used to define
\emph{Euclidean normal spaces}  at each point $P\in X_{\rm sm}$. Actually they considered a non-degenerate quadric in a hyperplane at $\infty$, essentially as we saw in the case of affine reciprocal polar varieties. Here we shall consider the orthogonality on all of $\mathbb P^n$, and we define the normal space at a smooth point $P$ as follows:
 \[N_PX:=\langle {T_P X}^\perp,P \rangle.\]
 
 We shall now see that, by passing to the Nash modification of  $X$, the Euclidean normal spaces are the fibers of a projective bundle.
 The Nash modification $\nu\colon\overline X\to X$ is the ``smallest'' proper, birational map such that the pullback of the cotangent sheaf $\Omega^1_X$ of $X$ admits a locally free quotient of rank $m$. This is equivalent to $\nu^*\mathcal P^1_X(1)$ admitting a locally free quotient of rank $m+1$. Denote this quotient by $\mathcal P$, and let $\mathcal K$ denote the kernel of the surjection $\mathcal O_{\overline X}^{n+1} \to \mathcal P$. Thus $\mathcal K$ is a modification of the conormal sheaf $\mathcal N_{X/\mathbb P^n}$ of $X$ in $\mathbb P^n$ twisted by $\mathcal O_X(1)$.

 The quadric $Q$ gives an isomorphism $\mathcal O_X^{n+1} \cong \mathcal (O_X^{n+1})^\vee$, hence we get a quotient $\mathcal O_{\overline X}^{n+1} \to \mathcal K^\vee $, whose (projective) fibers are the spaces $ T_P X^\perp$. Adding 
 the point map $\mathcal O_X^{n+1} \to \mathcal O_X(1)$, we get a surjection
 \[\mathcal O_{\overline X}^{n+1} \to \mathcal E:=\mathcal K^\vee \oplus \mathcal O_{\overline X}(1),\]
 whose (projective) fibers are the Euclidean normal spaces $N_PX$. Indeed, $\mathbb P(\mathcal E)\subset {\overline X}\times \mathbb P^n$, and the fibers of the structure map $\mathbb P(\mathcal E)\to {\overline X}$ above smooth points of $X$ are the spaces $N_PX\subset \mathbb P^n$ defined above. We call $\mathcal E$ the \emph{Euclidean normal bundle} of $X$ with respect to $Q$ (cf. \cite{C-T} and \cite{ED}).
 \medskip
 
  Let $p\colon \mathbb P(\mathcal E) \to {\overline X}$ denote the structure map, and let $q\colon \mathbb P(\mathcal E) \to \mathbb P^n$ denote the projection on the second factor. The map $q$ is  called the \emph{endpoint map} (for an explanation of the name, see  \cite{C-T}).
 \medskip

 Let $B\in \mathbb P^n$ be a (general) point. Then 
 \[ p(q^{-1}(B)) = \{ P\in X \, |\, B\in \langle T_PX^\perp,P \rangle \}.\]
 Letting $L:=B^\perp$, so that $B=L^\perp$, we see that
  \[ p(q^{-1}(B)) = \{ P\in X \, |\, L^\perp\in \langle T_PX^\perp,P \rangle \} = W_L^\perp(X)\]
  is a reciprocal polar variety. In particular,
  the degree of $q$ is just the degree of the reciprocal polar variety.

  \begin{example}
   Assume $X\subset \mathbb P^2$ is a (general) plane curve of degree $d$.
  The \emph{reciprocal} polar variety is the intersection of the curve with its reciprocal polar, which has degree $d$, so $q$ has degree $d^2$. 
\end{example}

 In \cite{ED}  the degree of the endpoint map $q\colon \mathbb P(\mathcal E) \to \mathbb P^n$ was called the \emph{Euclidean distance degree} of $X$:
 \[ {\rm ED}\deg X=p_*c_1(\mathcal O_{\mathbb P(\mathcal E)}(1))^n=s_m(\mathcal E),\]
 where $m=\dim X$ and $s_m$ denotes the $m$th Segre class.  The reason for the name is the relationship to computing critical points for the distance function in the Euclidean setting. We refer to  \cite{ED} for many more details. In the case of curves and surfaces (and in a slightly different setting), the degree of $q$ is called the normal class in \cite{JP2}.

 Since $\mathcal E=\mathcal K^\vee \oplus \mathcal O_{\overline X}(1)$, we get  \[s(\mathcal E)=s(\mathcal K^\vee)s(\mathcal O_{\overline X}(1))=c(\mathcal P)c(\mathcal O_{\overline X}(-1))^{-1}.\]
 Therefore, 
 \[s_m(\mathcal E)=\sum_{i=0}^m c_{i}(\mathcal P)c_1(\mathcal O_{\overline X}(1))^{m-i}.\]
 Since $c_{i}(\mathcal P)c_1(\mathcal O_{\overline X}(1))^{m-i}$ is the degree $\mu_i$ of the $i$th polar variety $[M_i]$ of $X$ \cite[p.~256]{Piene1}, we conclude (cf. \cite[Thm.~5.4]{ED}):
 \[  {\rm ED}\deg X=\sum_{i=0}^m \mu_i.\]
 The $\mu_i$ are called the \emph{ranks} (or classes) of $X$. Note that $\mu_0$ is the degree of $X$ and $\mu_{n-1}$ is the degree of the dual variety $X^\vee$ (provided the dimension of $X^\vee$ is $n-1$). It is known (see \cite[Prop. 3.6, p.~266]{Piene1}, \cite[3.3]{U}, and \cite[(4), p.~189]{K}) that the $i$th rank of $X$ is equal to the $(n-1-i)$th rank of the dual variety $X^\vee$ of $X$. As observed in \cite[Thm.~5.4]{ED}, it follows that the Euclidean distance degree of $X$ is equal to that of $X^\vee$.
 Moreover, whenever we have formulas for the ranks $\mu_i$, we thus get formulas for ${\rm ED}\deg X$.
 
 \begin{example}
 If $X\subset \mathbb P^n$ is a smooth hypersurface of degree $\mu_0=d$, then $\mu_i=d(d-1)^i$, hence in this case (cf. \cite[(7.1)]{ED})
 \[{\rm ED}\deg X=d\sum_{i=0}^{n-1} (d-1)^i= \frac{d((d-1)^n-1)}{d-2}.
 \]
 If $X$ has only isolated singularities, then only $\mu_{n-1}$ is affected, and we get (from Teissier's formula \cite[Cor. 1.5, p.~320]{Teissier}, and the Pl\"ucker formula for hypersurfaces with isolated singularities (\cite[II.3, p.~46]{Teissier2} and \cite[Cor. 4.2.1, p.~60]{Laumon})
  \[{\rm ED}\deg X= \frac{d((d-1)^n-1)}{d-2}-\sum_{P\in {\rm Sing}(X)}(\mu_P^{(n)}+\mu_P^{(n-1)}),
 \]
 where $\mu_P^{(n)}$ is the Milnor number and $\mu_P^{(n-1)}$ is the sectional Milnor number of $X$ at $P$.
 
 \end{example}
 
 \begin{example}
 Assume $X\subset \mathbb P^3$ is a generic projection of a smooth surface of degree $\mu_0=d$, so that $X$ has \emph{ordinary} singularities: a double curve of degree $\epsilon$, $t$ triple points, and $\nu_2$ pinch points. Then (using the formulas for $\mu_1$ and $\mu_2$ given in \cite[p.~18]{Piene2})
 \[{\rm ED}\deg X=\mu_0+\mu_1+\mu_2=d^3-d^2+d-(3d -2)\epsilon - 3t-2\nu_2.\]
 Further examples can be deduced from results in \cite{Piene1}, \cite{Piene2}, and \cite{Piene3}.
 \end{example}
 \medskip
  
  \section{Focal loci}
   
 The \emph{focal locus} (see e.g. \cite{C-T} for an explanation of the name,  or \cite{ED}, where it is denoted the \emph{ED discriminant})  is the branch locus  $\Sigma_X$ of the endpoint map $q\colon \mathbb P(\mathcal E)\to \mathbb P^n$. More precisely,  let $R_X$ denote the ramification locus of $q$; by definition, $R_X$ is the subscheme of  $\mathbb P(\mathcal E)$ given on the smooth locus $\mathbb P(\mathcal E)_{\rm sm}$ by the $0$th Fitting  ideal $F^0(\Omega^1_{\mathbb P(\mathcal E)/{\mathbb P^n}|_{\mathbb P(\mathcal E)_{\rm sm}}})$.  The focal locus $\Sigma_X$ is the closure of the image $q(R_X)$.
  
 Recall that we have, on the Nash modification $\nu\colon  \overline X\to X$, the exact sequence
 \[0\to \mathcal K\to \mathcal O_{\overline X}^{n+1} \to \mathcal P\to 0,\]
 where $\mathcal K$ and $\mathcal P$ are the Nash bundles of the sheaves $\mathcal N_{X/\mathbb P^n}(1)$ and $\mathcal P_X^1(1)$ respectively, and that $\mathcal E=\mathcal K^\vee \oplus \mathcal O_{\overline X}(1)$.
 
 Let $Z\to \overline X$ be a resolution of singularities, and, by abuse of notation, denote by $\mathcal K$, $\mathcal P$, $\mathcal E$ also their pullbacks to $Z$. The class $[R_X]$ of the closure of the ramification locus of $q\colon \mathbb P(\mathcal E) \to \mathbb P^n$ is given by
  \[[R_X]=c_1(\Omega^1_{\mathbb P(\mathcal E)})-q^*c_1(\Omega^1_{\mathbb P^n})=c_1(\Omega^1_{\mathbb P(\mathcal E)})+(n+1)c_1(\mathcal O_{\mathbb P(\mathcal E)}(1)).\]
  
 Using the exact sequences 
 \[0\to p^*\Omega^1_Z\to \Omega^1_{\mathbb P(\mathcal E)}Ê\to  \Omega^1_{\mathbb P(\mathcal E)/Z} \to 0\]
 and 
 \[0\to \Omega^1_{\mathbb P(\mathcal E)/Z} \to p^*\mathcal E\otimes O_{\mathbb P(\mathcal E)}(-1) \to \mathcal O_{\mathbb P(\mathcal E)} \to 0 \]
 we find
 \[[R_X] =p^*\bigl(c_1(\Omega_Z^1)+c_1(\mathcal P)+c_1(\mathcal O_Z(1))\bigr)+mc_1(\mathcal O_{\mathbb P(\mathcal E}(1)).\]
 Therefore the degree of $R_X$ with respect to the map $q$ is given by
 \[\deg R_X=\bigl(c_1(\Omega_Z^1)+c_1(\mathcal P)+c_1(\mathcal O_Z(1))\bigr)p_* c_1(\mathcal O_{\mathbb P(\mathcal E)}(1))^{n-1}+mp_*c_1(\mathcal O_{\mathbb P(\mathcal E)}(1))^{n}.\]
In terms of Segre classes of $\mathcal E$ this gives
\[\deg R_X=\bigl(c_1(\Omega_Z^1)+c_1(\mathcal P)+c_1(\mathcal O_Z(1))\bigr)s_{m-1}(\mathcal E)+ms_m(\mathcal E),\]
which gives, since $c_1(\Omega^1_Z)=c_1(\mathcal P^1_Z(1))-(m+1)c_1(\mathcal O_Z(1))$, 
\[\deg R_X=\bigl(c_1(\mathcal P^1_Z(1))+c_1(\mathcal P)-mc_1(\mathcal O_Z(1))\bigr)s_{m-1}(\mathcal E)+ms_m(\mathcal E).\]
Now $s_m(\mathcal E)=\sum_{i=0}^m\mu_i$ and  $s_{m-1}(\mathcal E)c_1(\mathcal O_Z(1))=\sum_{i=0}^{m-1}\mu_i$, hence
\[\deg R_X=\bigl(c_1(\mathcal P^1_Z(1))+c_1(\mathcal P)\bigr)s_{m-1}(\mathcal E)-m\sum_{i=0}^{m-1}\mu_i+m\sum_{i=0}^m\mu_i,\]
hence
\[\deg R_X=\bigl(c_1(\mathcal P_Z^1(1))+c_1(\mathcal P)\bigr)s_{m-1}(\mathcal E)++m\mu_m.\]

In the special case when $X\subset \mathbb P^n$ is a hypersurface ($m=n-1$), 
we know by \cite[Cor.~3.4]{Piene1} that
\[c_1(\mathcal P_Z^1(1))=c_1(\mathcal P)+c_1(\mathcal R_{Z/X}^{-1}),\]
where $\mathcal R_{Z/X}=F^0(\Omega^1_{Z/X})$ is the (invertible) ramification ideal of $Z\to X$. Hence we get
\[\deg R_X=\bigl(2c_1(\mathcal P)+c_1(\mathcal R_{Z/X}^{-1})\bigr)s_{n-2}(\mathcal E)+(n-1)\mu_{n-1},\]
or
\[\deg R_X=\bigl(2c_1(\mathcal P)+c_1(\mathcal R_{Z/X}^{-1})\bigr)\sum_{i=0}^{n-2}c_i(\mathcal P)c_1(\mathcal O_Z(1))^{n-2-i}+(n-1)\mu_{n-1}.\]

\begin{example}
Let $X\subset \mathbb P^2$ be a plane curve of degree $\mu_0$ and class $\mu_1=c_1(\mathcal P)$.  Then
\[ \deg R_X=2c_1(\mathcal P)+\kappa +\mu_1=3\mu_1+\kappa,\]
where $\kappa$ is the ``total number of cusps'' of $X$. Note that, again by \cite[Cor.~3.4]{Piene1}, $3\mu_1+\kappa = 3\mu_0+\iota$, where $\iota$ is the ``total number of inflection points'' of $X$. But the degree $\mu_0(X)$ of $X$ is equal to the class $\mu_1(X^\vee)$ of the dual curve $X^\vee$, and $\iota(X)$ of $X$ is $\kappa(X^\vee)$ of $X^\vee$. This shows that the degree of the focal locus, or ED discriminant, of the dual curve is equal to that of $X$:
\[ \deg R_{X^\vee}=3\mu_1(X^\vee)+\kappa(X^\vee)=3\mu_0+3\iota=\deg R_X.\]

The focal locus of a plane curve is also known as the \emph{evolute} or  the \emph{caustic by reflection}. So, provided the maps $R_X\to \Sigma_X$ and $R_{X^\vee}\to \Sigma_{X^\vee}$ are birational, we have shown that the degree of the evolute of $X$ is equal to the degree of the evolute of the dual curve $X^\vee$. For more on evolutes, see \cite{JP}. In the case that $X$ is a ``Pl\"ucker curve'' of degree $d=\mu_0$ and having only $\delta$ nodes and $\kappa$ ordinary cusps as singularities, and $\iota$ ordinary inflection points, then the classical formula, due to Salmon, is
 \[\deg R_X = 3d(d-1)-6\delta -8\kappa.\]
  Since in this case $\mu_1=d(d-1)-2\delta - 3\kappa$, this checks with our formula. Moreover, since  the number of inflection points is $\iota=3d(d-2)-6\delta - 8\kappa$, $\deg R_{X^\vee}=3d+\iota=3d(d-1)-6\delta -8\kappa =\deg R_X$, as it should.
\end{example}
\medskip

 If $X$ is smooth, we have $Z={\overline X}=X$, $\mathcal E=\mathcal N_{X/\mathbb P^n}(1)^\vee \oplus \mathcal O_X(1)$, and $s(\mathcal N_{X/\mathbb P^n}(1)^\vee)=c(\mathcal P^1_X(1))$. Hence we can compute  the class of $R_X$ in terms of the Chern classes of $X$ and $\mathcal O_X(1)$. 
We get
  \[\deg R_X=2\bigl(c_1(\Omega_X^1)\sum_{i=0}^{m-1}c_i(\mathcal P^1_X(1))c_1(\mathcal O_X(1))^{m-1-i}+(m+1)\sum_{i=0}^{m-1}\mu_i\bigr) +m\mu_m.\]
Since $c_i(\mathcal P^1_X(1))=\sum_{j=0}"\binom{m+1-i+j}{j}c_{i-j}(\Omega_X^1)c_1(\mathcal O_X(1))^j$, and since the $\mu_i$'s can be expressed in terms of the Chern numbers $c_{m-j}(\Omega_X^1)c_1(\mathcal O_X(1))^j$, we see that also $\deg R_X$ can be expressed in terms of these Chern numbers and the Chern numbers $c_1(\Omega_X^1)c_{m-1-j}(\Omega_X^1)c_1(\mathcal O_X(1))^j$.
  
\begin{example}
Assume $X\subset \mathbb P^n$ is  a smooth curve of degree $d$. Then
\[\deg R_X= 2(2g-2)+4\mu_0+\mu_1=2(2g-2)+4d+2d+2g-2=6(d+g-1),\]
as in \cite[Ex. 7.11]{ED}.
\end{example}

\begin{example}
Let $X\subset \mathbb P^n$ be a smooth surface of degree $\mu_0=d$. Then, as in \cite[Section 5]{C-T} we get:
\[\deg R_X=2(15d+9c_1(\Omega_X^1)c_1(\mathcal O_X(1))+c_1(\Omega_X^1)^2+c_2(\Omega_X^1)).\]
\end{example}

 \begin{example}
 Let $X\subset \mathbb P^n$ be a general hypersurface ($m=n-1$) of degree $\mu_0$. It is known that in this case $R_X\to \Sigma_X$ is birational \cite[Thm.~2]{T}.  
 Since $c_1(\Omega_X^1)=(\mu_0-n-1)c_1(\mathcal O_X(1))$ we get
  \[\deg \Sigma_X=\deg R_X=(2\mu_0-n-1)s_{n-2}(\mathcal E)c_1(\mathcal O_X(1))+(n-1)s_{n-1}(\mathcal E).\]
 Hence
 \[ \deg \Sigma_X=(2\mu_0-n-1)\sum_{i=0}^{n-2} \mu_i +(n-1)\sum_{i=0}^{n-1} \mu_i=(n-1)\mu_{n-1}+2(\mu_0-1)\sum_{i=0}^{n-2}\mu_i.\]
  For a smooth hypersurface of degree $d$ in $\mathbb P^n$, we have $\mu_i=d(d-1)^i$. Hence
  \[ \deg \Sigma_X=(n-1)d(d-1)^{n-1}+2d(d-1)((d-1)^{n-1}-1)(d-2)^{-1},\]
  which checks with the formula found in \cite[Thm.~2]{T}.
  \end{example}

\end{document}